\documentclass[10pt,conference]{IEEEtran}
\usepackage{graphicx}

\usepackage{cite}
\ifCLASSINFOpdf
\else
\fi

\usepackage{mathrsfs,amsfonts,amssymb,amsthm}
\usepackage[cmex10,intlimits]{amsmath}
\usepackage{array}

\usepackage{mdwmath}
\usepackage{mdwtab}

\begin{document}
\newtheorem{theorem}{Theorem}
\newtheorem{proposition}{Proposition}
\newtheorem{corollary}{Corollary}
\newtheorem{lemma}{Lemma}
\newtheorem{definition}{Definition}
\newtheorem{assumption}{Assumption}
\newtheorem{remark}{Remark}
\newtheorem{operation}{Operation}
\newcommand{\ams}{\textit{Ann. Math. Statist.}}
\newcommand{\dfn}{\stackrel{\triangle}{=}}
\newcommand{\argmax}{\operatornamewithlimits{arg\,max}}
\newcommand{\argmin}{\operatornamewithlimits{arg\,min}}
\newcommand{\argsup}{\operatornamewithlimits{arg\,sup}}
\newcommand{\arginf}{\operatornamewithlimits{arg\,inf}}

\title{Decentralized Two-Sided Sequential Tests for A Normal Mean}

\author{\IEEEauthorblockN{Yan Wang and Yajun Mei}
\authorblockA{School of Industrial and Systems Engineering\\
Georgia Institute of Technology,
Atlanta, Georgia 30332--0225\\
Email: \{ywang67, ymei\}@isye.gatech.edu}
}

\maketitle

\begin{abstract}
\boldmath
This article is concerned with decentralized sequential testing of a normal mean $\mu$ with two-sided alternatives. It is assumed that in a single-sensor network system with limited local memory, i.i.d. normal raw observations are observed at the local sensor, and quantized into binary messages that are sent to the fusion center, which makes a final decision between the null hypothesis $H_0: \mu = 0$ and the alternative hypothesis $H_1: \mu = \pm 1.$ We propose a decentralized sequential test  using the idea of tandem quantizers (or equivalently, a one-shot feedback). Surprisingly, our proposed test only uses the quantizers of the form $I(X_{n} \ge \lambda),$ but it is shown to be asymptotically Bayes. Moreover, by adopting the principle of  invariance, we also investigate decentralized invariant tests with the stationary quantizers of the form $I(|X_{n}| > \lambda),$ and show that $\lambda = 0.5$ only leads to a suboptimal decentralized invariant sequential test. Numerical simulations are conducted to support our arguments.
\end{abstract}

\IEEEpeerreviewmaketitle

\section{Introduction}\label{sec:intro}

Decentralized sequential detection has many important applications such as signal detection and sensor networks, see, for example, Blum, Kassam, and Poor \cite{blum97}.  Veeravalli, Basar, and Poor \cite{vee} characterizes the Bayesian solutions in the system with limited local memory and full feedback. Recently Mei \cite{mei} develops the first complete asymptotic theory for decentralized sequential detection. However,  existing research only focuses on the simplest model when both null and alternative hypotheses are completely specified.

In this article, we will consider a more flexible model of decentralized sequential detection  in which hypotheses are composite. To highlight our main ideas, we focus on the following specific problem in a single-sensor network system, since the extension to the system with multiple (conditionally independent) sensors is straightforward. Assume that the single local sensor observes a sequence of raw observations $X_{1}, X_{2}, \cdots$ over time $n$ and the $X_{n}$'s are i.i.d. having a normal distribution $N(\mu,1)$. Suppose we are interested in testing
\begin{equation}\label{equ:Hypothesis}
H_0: \mu=0\ \ \ \textrm{versus}\ \ \ H_1: \mu=\pm 1.
\end{equation}
In the centralized context, one uses the raw observations $X_{i}$'s to decide which of $H_0$ and $H_1$ is true, and such a problem has been well studied in the mature field of sequential analysis (Wald \cite{wald47}).
In the context of decentralized detection,  due to data compression and communication constraints, the local sensor needs to quantize the data $X_{n}$'s and send a binary message $U_{n} \in \{0, 1\}$ to the fusion center, which then utilize the quantized messages $U_{n}$'s to decide which of $H_0$ and $H_1$ is true. Following Veeravalli, Basar, and Poor \cite{vee} and Mei \cite{mei}, it is assumed that at time $n$, the quantized message $U_n$ sent from the sensor to the fusion center only depends on the current raw observation $X_{n}$ and possibly feedback from the fusion center. In other words, at time $n,$ quantized message $U_n$ satisfies
\begin{equation}\label{equ:lmtdlclMmry}
U_n=\phi_n(X_n; V_{n-1}) \qquad \in\{0,1\},
\end{equation}
where the feedback $V_{n-1}$ only depends on past sensor messages:
$V_{n-1}=\psi_{n}(U_{[1,n-1]});$ where $U_{[1,n-1]}=(U_1,\dots,U_{n-1})$.

In the decentralized sequential detection problems, one wants to determine how to design sensor quantizers $\{\phi_n\}$ in  (\ref{equ:lmtdlclMmry}) and how to make a sequential decision at the fusion center, so that the overall performance of the system is optimal (in some suitable senses). A central challenge is to determine the form of (binary) quantizers $\phi_n$'s for (asymptotically) optimal decentralized tests. In the simplest model when both null and alternative hypotheses are completely specified, the best quantizers are of the form of monotone likelihood ratio quantizers (MLRQ), see Tsitsiklis \cite{tsi}, and take the following simple form in the case of testing normal means:
\begin{equation}\label{equ:QntzrTypeI}
\phi(X)=I(X\ge \lambda).
\end{equation}
When the hypothesis is composite, the MLRQ is no longer applicable. In particular, it is unclear whether the quantizers in (\ref{equ:QntzrTypeI}) still lead to (asymptotically) optimal decentralized solutions when testing the hypotheses in (\ref{equ:Hypothesis}). Indeed, our intuition may suggest us that a (more) attractive candidate can be $\phi(X) = I( |X| \le 0.5),$ or more generally,
$\phi(X)=I(\lambda_1\le X\le \lambda_2).$
Moreover, it is unclear whether other more complicated quantizers are necessary or not.

In this article, we tackle the form of binary quantizers by using the concept of unambiguous likelihood quantizer (ULQ) proposed by Tsitsiklis \cite{tsi} (the MLRQ is a special case of ULQ). Surprisingly, by combining the ULQ with the idea of tandem quantizer in Mei \cite{mei}, we show that at most one switch between two different quantizers of the form in (\ref{equ:QntzrTypeI}) is sufficient to construct the asymptotically optimal decentralized sequential test when testing the composite hypotheses in (\ref{equ:Hypothesis}).

Observing the symmetries of the densities, it is also natural to adopt the principle of invariance (see, for example, Lehmann \cite{leh}). Specifically, if we consider $|X_n|,$ the problem of testing hypothesis in (\ref{equ:Hypothesis}) becomes one of testing a simple null versus a simple alternative hypothesis on $|X_n|$. This viewpoint allows us to apply the asymptotic theory in Mei \cite{mei} to investigate decentralized {\it invariant} sequential tests. It is interesting to note that among stationary quantizers of the form  $\phi(X) = I( |X| \ge \lambda),$ the intuitive choice of $\lambda = 0.5$ leads suboptimal decentralized  invariant sequential tests.

The remainder of this article is organized as follows. Section II provides a formal mathematical formulation of decentralized sequential hypothesis testing problem. In Section III, we propose a family of decentralized sequential tests, and proves its asymptotic optimality properties. The nontrivial part of the proof is in Subsection III.C, which characterize optimal quantizers via unambiguous likelihood quantizers (ULQ).  Section IV focuses on the decentralized {\it invariant} sequential tests. Section V reports numerical simulations.

\bigskip
\section{Problem Formulation}

Assume that the raw data $X_1, X_2, \cdots$ are i.i.d. with $N(\mu, 1),$ and suppose that over time $n,$ the quantized message $U_{n} \in \{0, 1\},$ defined in  (\ref{equ:lmtdlclMmry}), only depends on $X_{n}$ and possibly $V_{n-1},$ a feedback summarizing past history of all quantized messages $U_{[1,n-1]}.$  Here we intentionally do not put any restrictions on the range or frequency of the feedback $V_{n-1}$, as it turns out that a simple one-shot feedback is sufficient to construct an asymptotic optimal solution.

For three hypotheses in (\ref{equ:Hypothesis}), denote by $f$, $g_1$ and $g_2$ the probability densities of $X_{n}$'s when $\mu=0,-1,1$ respectively. Also denote the corresponding probability measures and expectations by $\{\textbf{P}_{f}, \textbf{E}_f\}$, $\{\textbf{P}_{g_1},\textbf{E}_{g_1}\}$ and $\{\textbf{P}_{g_2},\textbf{E}_{g_2}\}$. Assume a priori distribution $\mathbf{\pi}=(\pi_f,\pi_{g_1}, \pi_{g_2})$ is assigned to the three states of nature, and let
\[
\textbf{P}_\pi=\pi_f \textbf{P}_f+\sum_{i=1}^{2}\pi_{g_i}\textbf{P}_{g_i};\ \ \textbf{E}_\pi=\pi_f \textbf{E}_f+\sum_{i=1}^{2}\pi_{g_i}\textbf{E}_{g_i}.
\]

To characterize a decentralized sequential test $\delta$, denote by $N$ the time when the test $\delta$ decides to stop taking observations, i.e., $N$ is the sample size of $\delta$. Once stopped, the test (the fusion center) makes a decision $d\in\{0,1\}$, corresponding to $H_0$ and $H_1$, based on the information it receives up to that time. In summary, a decentralized sequential test  $\delta$ includes a sequence of quantizers $\{\phi_1,\phi_2,\dots\}$, a sequence of feedback functions $\{\psi_1,\psi_2,\dots\}$, a stopping time $N$ at the fusion center, and a decision function $d \in\{0,1\}$.

As in Wald \cite{wald47} and Veeravalli, Basar, and Poor \cite{vee}, define a Bayes risk of
a decentralized sequential test $\delta$  as
\begin{multline}\label{equ:ttlCstsApriorPi}
\mathcal{R}_c(\delta)=\pi_f[c\textbf{E}_f(N)+W_f\textbf{P}_f\{d=1\}]\\
	+\sum_{i=1}^{2}\pi_{g_i}[c\textbf{E}_{g_i}(N)+W_{g_i}\textbf{P}_{g_i}\{d=0\}],
\end{multline}
with $c$ the incremental cost of each sample and $\{W_f, W_{g_1}, W_{g_2}\}$ cost of making incorrect decisions. The Bayes formulation of decentralized sequential hypothesis testing problems can then be stated as follows.

\medskip
{\it Problem (P1):}  Minimize the Bayes risk $\mathcal{R}_c(\delta)$ in (\ref{equ:ttlCstsApriorPi}) among all possible decentralized sequential tests.
\medskip

Let $\delta_{B}^*(c)$ denote a Bayes solution to the decentralized sequential detection problem, i.e., $\delta_{B}^*(c)=\argmin_{\delta}\{\mathcal{R}_c(\delta)\}.$ Since it is extremely difficult, if possible at all, to find the exact form of $\delta_B^*(c)$ when hypotheses are composite, we adopt the asymptotically optimal approach, i.e., to find a family of decentralized tests $\{\delta(c)\}$ such that
\begin{equation*}\label{equ:dfnAspBays}
\lim_{c \rightarrow 0} \mathcal{R}_c(\delta_{B}^*(c))/\mathcal{R}_c(\delta(c)) =  1.
\end{equation*}

\medskip
\section{Our Proposed Test $\delta_{I}(c)$}\label{sec:ABD}

In this section we propose a family of tests $\{\delta_I(c)\}$ that is asymptotically Bayes. Our proposed test is a two-stage procedure, and it assumes that the fusion center will send a one-shot feedback $V$ taking values in $\{0,1,2\}$, representing a preliminary decision on $f, g_1$ or $g_2.$

\subsection{Definition of Test $\delta_{I}(c)$.}

Our proposed test $\delta_{I}(c)$ is defined as follows.

1) \textit{First Stage:} Choose positive values $u(c)<1/2$ satisfying
\begin{equation}\label{equ:stage1thresh}
u(c)\to 0 \  \quad \mbox{and} \quad \frac{\log u(c)}{\log c}\to 0.
\end{equation}
In the first stage, the local sensor quantizes the raw data $X_{n}$'s by a stationary quantizer $\phi^0(X)=I(X\ge 0)$. Based on the quantized message $U_{n}=\phi^0(X_{n})$ at each time $n,$ the fusion center updates the posterior distribution of the three densities
$(\pi_{f,n},\pi_{g_1,n}, \pi_{g_2,n})$ recursively. For example when $U_n=1$,
	  \begin{equation*}\label{equ:updtsaPos}
\pi_{f,n}=\frac{\pi_{f,n-1}\textbf{P}_f(U_n=1)}{\sum_{i\in\{f,g_1,g_2\}}\textbf{P}_i(U_n=1)}.
	  \end{equation*}
The fusion center stops the first stage at time
\begin{equation*}\label{equ:stage1stoppingtime}
N_1=\min\{n \ge 1: \max\{\pi_{f,n},\pi_{g_1,n},\pi_{g_2,n}\}\ge 1-u(c)\},
\end{equation*}
and makes a preliminary decision $d^0\in\{f,g_1,g_2\}$ at time $N_1$ satisfying
\begin{equation*}\label{equ:stage1decision}
\pi_{d^0,N_1}=\max\{\pi_{f,N_1},\pi_{g_1,N_1},\pi_{g_2,N_1}\}.
\end{equation*}

2) \textit{Second Stage:} In this stage, it is essential for the sensor to switch to one of the following three ``optimal" quantizers, depending on the preliminary decision of the fusion center in the first stage:
\begin{eqnarray*}
\phi^*_f(X) &=& I(X\ge 0);\\
\phi^*_{g_1}(X) &=& I(X\ge -0.7941);\\
\phi^*_{g_2}(X) &=& I(X\ge 0.7941).
\end{eqnarray*}

Specifically, after the fusion center stops at the first stage, it will send its preliminary decision $d^{0}\in\{f,g_1,g_2\}$ back to the local sensor as a one-shot feedback $V$. Then the local sensor switch to the above optimal quantizer $\phi_{d^0}^*$ and use it for all incoming raw data.

In the second stage, the fusion center continues to recursively update the posterior distributions with data $U_{N_1+1},U_{N_1+2},\dots$ and it decides to stop the second stage at time
\begin{equation*}
N=\min\{n\ge N_1: \dfrac{\pi_{f,n}W_f}{\sum_{i=1}^{2}\pi_{g_i,n}W_{g_i} }\not\in(c,\frac{1}{c})\}\label{equ:stoppingtime}
\end{equation*}
with a final decision $d=0$ ($d=1$) if upper (lower) bound is crossed.

\subsection{Asymptotic Optimality of $\delta_{I}(c)$.}\label{sec:AOP}

The asymptotic optimality properties of test $\delta_{I}(c)$ are summarized in the following theorem and its corollary:
\begin{theorem}\label{the:lowerBdReached}
For any decentralized sequential tests $\{\delta(c)\},$ if
\begin{equation}\label{equ:lowerBdfCond}
\textbf{P}_{\varphi}(\textrm{decision incorrect})=O(c\log c),
\end{equation}
for $\varphi\in\{f,g_1,g_2\},$ then the stopping time $N$ of $\delta(c)$ satisfies
\begin{equation}\label{equ:lowerBdf}
\begin{gathered}
\textbf{E}_f(N) \ge (1+o(1))|\log c|/0.3137\\
\textbf{E}_{g_1}(N)\ge (1+o(1))|\log c|/0.3186\\
\textbf{E}_{g_2}(N)\ge (1+o(1))|\log c|/0.3186,
\end{gathered}
\end{equation}
and our proposed tests $\{\delta_{I}(c)\}$ attain all three lower bounds in (\ref{equ:lowerBdf}) simultaneously.
\end{theorem}

\begin{corollary}\label{cor:BayesOptDeltBc}
Tests $\{\delta_{I}(c)\}$ are asymptotically Bayes. Moreover, both  $\{\delta_{I}(c)\}$  and Bayes solution $\delta_{B}^*(c)$ satisfy as $c\to 0$,
\begin{multline*}\label{equ:BayesOptimal}
\mathcal{R}_c(\delta_{I}(c))=c|\log c|(1+o(1)) \Big[\frac{\pi_f}{0.3137} + \frac{\pi_{g_1}+\pi_{g_2}}{0.3186}\Big].
    \end{multline*}

\end{corollary}

Note that the asymptotic optimality properties in Theorem \ref{the:lowerBdReached} do not depend on either the priori distribution $\pi$ or the loss for incorrect decisions $\{W_f, W_{g_1}, W_{g_2}\}$. This is consistent with the centralized sequential hypothesis testing, see, Chernoff \cite{che}.

Before proving Theorem \ref{the:lowerBdReached}, let us first introduce some necessary notation. Denote by $\Phi$ the set of deterministic quantizers that consists of all (deterministic) measurable functions from $\mathbb{R}$ to $\{0,1\}$. Define a ``random quantizer'' $\bar{\phi}$ as a probability measure that assigns certain masses $\{p_i\}$ on a finite subset $\{\phi_i\}\in\Phi$, and denote by $\bar{\Phi}$ the set of all quantizers, deterministic or random. Note that a deterministic quantizer can be thought of as a special case of random quantizer that assigns a probability of $1$ to itself.

In the context of decentralized detection, we adopt the following implementation for a random quantizer $\bar{\phi}:$ The fusion center first selects a deterministic quantizer $\phi\in\Phi$ randomly according to the probability measure $\{p_i\}$ assigned by $\bar{\phi},$ and then the local sensor quantizes the raw data by the chosen deterministic quantizer $\phi.$ Such a procedure repeats whenever the local sensor uses $\bar{\phi}$ to quantize a new raw observation. We want to emphasize that it is essential to assume that if a random quantizer $\bar{\phi}$ is applied, the fusion center retains the information about which deterministic quantizer it chooses (otherwise the fusion center will lose significant information).

Observe that for a given deterministic quantizer $\phi$, the K-L information number $I^{\phi}(f,g_1)$ is
\begin{equation*}\label{equ:dfnPhiKLInfo}
I^{\phi}(f,g_1)=\sum_{i=0}^{1}\textbf{P}_f(\phi(X)=i)\log \frac{\textbf{P}_f(\phi(X)=i)}{\textbf{P}_{g_1}(\phi(X)=i)}.
\end{equation*}
With our implementation of random quantizers, for a given random quantizer $\bar{\phi}$ that assigns probability mass $p_1,\dots,p_n$ onto $\phi_1,\dots,\phi_n$, it is easy to see that the corresponding K-L information number $I^{\bar{\phi}}(f,g_1)$ (at the fusion center) is
\begin{equation}\label{equ:DfnKLInfoNORdmQntz}
I^{\bar{\phi}}(f,g_1)=\sum_{i=1}^{n}p_i I^{\phi_i}(f,g_1).
\end{equation}
Similarly, we can also define the quantities $I^{\phi}(f,g_2)$, $I^{\bar \phi}(f,g_2)$, or $I^{\phi}(g_{i},f)$, $I^{\bar \phi}(g_{i},f)$ for $i=1,2$.

We are now ready to rigorously define the optimal quantizers and the corresponding K-L information number. Define the optimal quantizer with respect to $g_{i}$ as
\[\bar{\phi}_{g_i}=\argsup_{\bar{\phi}\in{\bar{\Phi}}}\{I^{\bar{\phi}}(g_i,f)\}, \ \ i=1,2\]
and define the optimal quantizer with respect to $f$ as
\begin{equation}\label{equ:OptPhif}
\bar{\phi}_f=\argsup_{\bar{\phi}\in\bar{\Phi}}\{\min\{I^{\bar{\phi}}(f,g_1),I^{\bar{\phi}}(f,g_2)\}\}.
\end{equation}

Moreover, define the corresponding K-L information number of these two quantizers as $I_{g_i}=I^{\bar{\phi}_{g_i}}(g_i,f)$, $i=1,2$ and $I_f=\min\{I^{\bar{\phi}_f}(f,g_1),I^{\bar{\phi}_f}(f,g_2)\}$.

With these notation, let us state the following proposition without proof, as it is just a special case of Theorem 2 of Chernoff \cite{che} and Section V of Kiefer and Sacks \cite{ks}.

\begin{proposition}\label{prop:lowerBdReached}
For decentralized sequential tests  $\{\delta(c)\}$ satisfying (\ref{equ:lowerBdfCond}) in Theorem \ref{the:lowerBdReached}, for $\varphi\in\{f,g_1,g_2\}$, as $c\to 0$,
\begin{equation}\label{equ:proplowerBdf}
\textbf{E}_\varphi(N)\le (1+o(1))|\log c|/I_\varphi.
\end{equation}
To achieve the lower bounds in (\ref{equ:proplowerBdf}) simultaneously, one only needs to use the two-stage procedure as described in Section \ref{sec:ABD} for test $\delta_{I}(c)$, but for the second stage, $\phi^*_f$, $\phi^*_{g_1}$ and $\phi^*_{g_2}$ should be substituted by the optimal (random) quantizers $\bar{\phi}_f$, $\bar{\phi}_{g_1}$ and $\bar{\phi}_{g_2}$ respectively.
\end{proposition}

A comparison of Theorem \ref{the:lowerBdReached} and Proposition \ref{prop:lowerBdReached} shows that to prove Theorem \ref{the:lowerBdReached}, it suffices to show that in the context of testing a normal mean stated in (\ref{equ:Hypothesis}), the optimal quantizers $\bar{\phi}_f, \bar{\phi}_{g_1}, \bar{\phi}_{g_2}$ become
$\phi^*_f,  \phi^*_{g_1}, \phi^*_{g_2}$ described in Section \ref{sec:ABD}. Since $\phi^*_{g_1}$ or $\phi^*_{g_2}$ only involves two densities, the corresponding result follows immediately from the optimality of MLRQ's established in Tsitsiklis \cite{tsi}. Therefore, it remains to show that
\begin{eqnarray} \label{equ:barphifIsphiStarf}
\bar{\phi}_f=\phi^*_f,
\end{eqnarray}
which will be proved in the next subsection.

\subsection{Optimal Quantizer with respect to $f$}

The main objective of this subsection is to prove (\ref{equ:barphifIsphiStarf}), i.e., the optimal  (randomized) quantizer $\bar{\phi}_f(x)$ with respect to $f$ becomes the deterministic quantizer $\phi^*_f(x) = I(x \ge 0)$ when $f= N(0,1), g_1=N(-1,1)$ and $g=N(1,1).$

To prove this, for a given deterministic quantizer $\phi\in\Phi$, define for $i=0,1$,
\[q_i(\phi|\varphi)=\textbf{P}_\varphi(\phi(X)=i) \textrm{ and } q(\phi|\varphi)=(q_0(\phi|\varphi),q_1(\phi|\varphi)),\]
where $\varphi\in\{f,g_1,g_2\}$, and denote
\begin{equation}\label{equ:dfn_q_phi}
q(\phi)=(q(\phi|f); q(\phi|g_1); q(\phi|g_2)),
\end{equation}
then $q(\phi)$ completely characterizes the distribution of quantized message induced by the deterministic quantizer $\phi$, in the sense that if  $q(\phi_1)=q(\phi_2)$, the quantized data $\phi_1(X)$ and $\phi_2(X)$ have the same distribution, which implies that $I^{\phi_1}(f,g_1)=I^{\phi_2}(f,g_2).$

Let
\begin{equation*}\label{equ:QntzerVcts}
Q=\{q(\phi),\phi\in\Phi\}
\end{equation*}
be a subset of $\mathbb{R}^6$ and $\tilde{Q}$ be the convex hull of $Q$. (Here we do not use the usual symbol $\bar{Q}$ to avoid confusion with $\bar{\Phi}$.) For $\tilde{q}\in\tilde{Q}$, as in (\ref{equ:DfnKLInfoNORdmQntz}) and (\ref{equ:dfn_q_phi}), define $\tilde{q}=(\tilde{q}(f); \tilde{q}(g_1); \tilde{q}(g_2))$ and
\begin{equation}\label{equ:KLInfo_Qtilde}
I^{\tilde{q}}(f,g_1)=\sum_{i=0}^{1}\tilde{q}_i(f)\log\frac{\tilde{q}_i(f)}{\tilde{q}_i(g_1)},
\end{equation}
and quantities such as $I^{\tilde{q}}(f,g_{i})$ with $i=1$ or $2$ in an obvious extension. Note that the K-L definition in (\ref{equ:KLInfo_Qtilde}) is consistent with that in (\ref{equ:DfnKLInfoNORdmQntz}), since for a deterministic quantizer $\phi,$ we have $I^{q(\phi)}(f,g_1)=I^{\phi}(f,g_1)$.
	
Now let us state the concept of unambiguous likelihood quantizer (ULQ) proposed in Tsitsiklis \cite{tsi}. Let
\[v_i(X)=\frac{g_i(X)/f(X)}{1+g_1(X)/f(X)+g_2(X)/f(X)},\ \ i=1,2.\]
In our context (also see Lemma \ref{lem:ExtQGivelargerInfoNo} below), a  quantizer $\phi\in\Phi$ is a ULQ if there exists real number $a_0, a_1, a_2$ such that
\begin{equation}\label{equ:dfnULQ}
\phi(X)=I(a_0+a_1v_1(X)+a_2v_2(X)>0)
\end{equation}
and for $\varphi\in\{f,g_1,g_2\}$,
\begin{equation}\label{equ:ULQ0probStaysOnTheBoerder}
\textbf{P}_\varphi(a_0+a_1v_1(X)+a_2v_2(X)=0)=0.
\end{equation}
\begin{lemma}\label{lem:ULQCharacter}
Let $\phi\in\Phi$ be a ULQ, then up to a permutation of the values it takes, w.p.$1$ (under all $f, g_1, g_2$):
\begin{equation}\label{equ:ULQcharacter}
\phi(X)=I(\lambda_1\le X\le \lambda_2)\ \textrm{ or }\ \phi(X)=I(X\ge\lambda).
\end{equation}
\end{lemma}

\begin{proof}
As $X$ goes from $-\infty$ to $\infty$, both $1-v_1$ and $v_2$ strictly increase from $0$ to $1$. Hence it suffices to show $\phi(X)=I(t_1\le v_1(X)\le t_2)$ with $0\le t_1\le t_2\le 1$.

By (\ref{equ:dfnULQ}) and (\ref{equ:ULQ0probStaysOnTheBoerder}), ULQ's can be interpreted as: draw a straight line on $\mathbb{R}^2$ which intersects $\{(v_1(X),v_2(X))\}$ at a zero-probability set (under all $f, g_1, g_2$), $\phi$ will take value $0$ if $(v_1,v_2)$ stays in one side of the line and take value $1$ if it stays in the other. Since $\frac{d^2 v_2}{dv_1^2}<0$ for any $0<v_1<1$, a line intersects $\{(v_1(X),v_2(X))\}$ at at most two points, and thus relation (\ref{equ:ULQcharacter}) holds.
\end{proof}

The following lemma shows that the best quantizers can be found from the class of the ULQ's.
\begin{lemma}\label{lem:ExtQGivelargerInfoNo}
For a given quantizer $\bar{\phi}\in\bar{\Phi}$, there always exists another quantizer $\bar{\phi}' \in\bar{\Phi}$ that assigns probability masses only to the ULQ's and
\begin{equation}\label{equ:ExtremalPtsCombineGood}
I^{\bar{\phi}'}(f,g_i)\ge I^{\bar{\phi}}(f,g_i); \ \ i=1,2.
\end{equation}
\end{lemma}

\begin{proof}
\noindent Let $Q_\alpha$ be the exposed points of $\tilde{Q}$, then by Corollary 5.1 of Tsitsiklis \cite{tsi}, $q\in Q_\alpha$ if and only if there exists a ULQ $\phi$ such that  $q=q(\phi)$. By the compactness of the set $\{(v_1(X),v_2(X))\},$ it is straightforward to show that $Q_\alpha$ is identical with the extremal points of $\tilde{Q}$.  

From  (\ref{equ:DfnKLInfoNORdmQntz}), it is sufficient to prove (\ref{equ:ExtremalPtsCombineGood}) for deterministic $\phi$.  By the extreme properties of the ULQ's, there exist ULQ's $\phi_1,\dots,\phi_n$ and positive number $p_1,\dots,p_n$ such that
\[q(\phi)=\sum_{k=1}^{n}p_kq(\phi_k); \ \ \sum_{k=1}^{n}p_k=1.\]
Let $\bar{\phi}'$ be a random quantizer assigning mass $p_i$ to $\phi_i$ for $i=1,\dots,n$. By (\ref{equ:DfnKLInfoNORdmQntz}), for $i=1,2$,
	\[I^{\bar{\phi}'}(f,g_i)=\sum_{k=1}^{n}p_k I^{\phi_k}(f,g_i)\ge I^{q(\phi)}(f,g_i)=I^{\phi}(f,g_i).\]
\end{proof}

\begin{proof}[Proof of Relation (\ref{equ:barphifIsphiStarf})]
By Lemma \ref{lem:ExtQGivelargerInfoNo}, $\bar{\phi}_f$ can be achieved by randomizing ULQ's of the form $\phi(X)=I(X\ge \lambda)$ or $\phi(X)=I(\lambda_1\le X\le \lambda_2)$. Since $\sup_{\bar{\phi}\in\bar{\Phi}}\{\min\{I^{\bar{\phi}}(f,g_1),I^{\bar{\phi}}(f,g_2)\}\}$ must be reached on the boundary of a two dimensional convex set, it suffices to focus on quantizers that randomize between at most two ULQ's.

Numerically, we can simply optimize over discrete (finite) sets $\Pi_1=\{\lambda_i: \lambda_{i} \in \mathbb{R} \}$ and $\Pi_2=\{(\lambda_{1,i},\lambda_{2,i}):\lambda_{1,i}<\lambda_{2,i}\}.$ For each quantizer $\bar{\phi}$ that randomize between at most two ULQ's (with the values of $\lambda$ and $(\lambda_1,\lambda_2)$ in $\Pi_1$ or $\Pi_2,$ respectively), we can calculate the value 
\[\min\{I^{\bar{\phi}}(f,g_1),I^{\bar{\phi}}(f,g_2)\},\] 
and the quantizer $\bar{\phi}$ with the maximum value will be the approximation of the best quantizer. Our numerical computations support the optimality of quantizer $\phi^*_f(X)=I(X>0)$ in the sense of (\ref{equ:OptPhif}) to the precision of four digits after the decimal point.
\end{proof}

\section{Invariant Tests}\label{sec:alter}

One popular approach to tackle hypothesis testing with composite hypotheses is the principle of invariance, see for example Lehmann \cite{leh}. In our case, the two densities in $H_1$ are  reflections to each other, so if we pretend that $\{|X_n|\}$ are taken as the raw data, the problem in (\ref{equ:Hypothesis}) becomes a simple hypothesis-testing problem with
\begin{equation*}\label{equ:HypothesisInvariantTests}
H_0: \ |X_n|\sim \tilde{f} \quad \mbox{and} \quad H_1: |X_n|\sim \tilde{g},
\end{equation*}
where $\tilde{f}$ and $\tilde{g}$ are probability densities of the forms:
\begin{eqnarray*}\tilde{f}(x)&=& \frac{2}{\sqrt{2\pi}}e^{-\frac{x^2}{2}} 1\{x\ge 0\} \label{equ:absDensitiesF}
\\
\tilde{g}(x)&=& \frac{1}{\sqrt{2\pi}}(e^{-\frac{(x-1)^2}{2}}+e^{-\frac{(x+1)^2}{2}})1\{x\ge 0\}.\label{equ:absDensitiesG}
\end{eqnarray*}
Therefore, with $|X_n|$, we can develop ``good" decentralized invariant sequential tests based on the asymptotic optimality theory in Mei \cite{mei}. Below we will use the same notation as in previous sections, for instance, denote by $N$ the sample size, and denote by $\pi_{\tilde{f},n}$ and $\pi_{\tilde{g},n}$ the posterior distributions, resp. Similarly, $\{W_{\tilde{f}}, W_{\tilde{g}}\}$ is the cost of making incorrect decisions.

Let us consider the decentralized invariant test with {\it stationary} quantizers of the form $U_{n} = I(|X_{n}|\le\lambda)$. In this case, the quantized sensor messages $U_{n}$'s are i.i.d. (conditioned on each hypothesis) and the fusion center faces a classical sequential hypothesis testing problem. Thus the optimal policy at the fusion center is an SPRT based on $U_{n}$'s. That is, the  fusion center stops taking observations at time
\begin{equation*}
N=\min\{n \ge 1: \frac{\pi_{\tilde{f},n}W_f}{ \pi_{\tilde{g},n}W_{\tilde{g}}} \not\in (c, \frac{1}{c})\},
\label{equ:stoppingtimeAbsXSta}
\end{equation*}
and decides $H_0$ ($H_1$) is true if upper (lower) bound is crossed. Hence, in the following we will we pay special attention on how to choose a quantizer of the form $U_{n} = I(|X_{n}|\le\lambda).$

Under our setting, one natural choice of quantizer is
\begin{equation}\label{equ:StatQntzr}
\phi(|X|)=I(|X|\le 0.5)
\end{equation}
and denote by $\{\delta_{II}(c)\}$  the corresponding decentralized test with an SPRT at the fusion center.

By asymptotic optimality theory in Mei \cite{mei}, a better choice of $\lambda$  is to find one value that minimizes $$\frac{\pi_{\tilde f}}{I(\tilde f_{\lambda}, \tilde g_{\lambda})} + \frac{\pi_{\tilde g}}{I(\tilde g_{\lambda}, \tilde f_{\lambda})},$$ where $\tilde f_{\lambda}$ and $\tilde g_{\lambda}$ are the probability mass functions induced on $U_{n} = I(|X_{n}|\le\lambda)$ when the distribution of $|X_{n}|$ is $\tilde f$ and $\tilde g$. A simple numerical simulation shows that the best value is $\lambda = 1.2824,$ and
\begin{equation}\label{equ:BstStatQntzr}
\phi^*(|X|)=I(|X|\le 1.2824).
\end{equation}
Denote by $\{\delta_{III}(c)\}$  the corresponding decentralized test with an SPRT at the fusion center.

\section{Simulation}\label{sec:simul}

In this section we compare the three tests proposed in previous sections through numerical simulation.  We fix a priori distribution $\pi_{f}=\pi_{g_1}=\pi_{g_2}=1/3$ in our problem (by Theorem \ref{the:lowerBdReached}, this is not essential). This leads to $\pi_{\tilde f} = 1/3$ and $\pi_{\tilde g} = 2/3$ for the invariant tests in Section \ref{sec:alter}. In our simulations, the cost of making incorrect decisions are assumed to be $1$, and we consider three different values for the cost of taking an observation: $c=10^{-2}, 10^{-3}, 10^{-4}.$ In our proposed test  $\delta_{I}(c),$ it has an additional parameter $u(c)$ satisfying the conditions in (\ref{equ:stage1thresh}), and in our simulations we assume that $u(c) = 0.1.$

Table \ref{tab:sampleSize} and \ref{tab:error} report numerical simulations on $\textbf{P}(\textrm{decision incorrect})$ and the expected sample sizes $\textbf{E}(N).$  Since the probabilities of incorrect decisions are small, we use the importance sampling approach to simulate $\textbf{P}(\textrm{DI})$.

\begin{table}[!h]
\renewcommand{\arraystretch}{1.3}
\caption{Expected Values of Sample Sizes $\textbf{E}(N)= \pi_{f} E_{f}(N) + \pi_{g_1} E_{g_1}(N) + \pi_{g_2} E_{g_2}(N)$}
\label{tab:sampleSize}
\centering
\begin{tabular}{c||c|c|c}
\hline
 $\textbf{E}(N)$ & $c=10^{-2}$ & $c=10^{-3}$ & $c=10^{-4}$ \\
\hline
\hline
$\delta_{I}(c)$ & $20.2\pm0.2$ & $28.4\pm0.2$ & $36.3\pm0.2$ \\
	 \hline
$\delta_{II}(c)$ &$94.1\pm0.7$&$146.0\pm1.0$& $196.4\pm1.0$\\
	 \hline
$\delta_{III}(c)$ & $45.7\pm0.5$ & $69.0\pm0.5$ & $92.2\pm0.5$\\
	 \hline
\end{tabular}

\bigskip
\caption{Probabilities Of Making Incorrect Decisions}
\label{tab:error}
\centering
\begin{tabular}{c||c|c|c}
\hline
 $\textbf{P}(\textrm{DI})$ & $c=10^{-2}$ & $c=10^{-3}$ & $c=10^{-4}$ \\
\hline
\hline
$\delta_{I}(c)$ & $4.58\pm0.03$e-3 & $4.42\pm0.03$e-4 & $4.61\pm0.03$e-5 \\
	 \hline
$\delta_{II}(c)$ &$8.84\pm0.02$e-3&$8.85\pm0.02$e-4&$8.84\pm0.02$e-5\\
\hline
$\delta_{III}(c)$ & $8.02\pm0.02$e-3 & $8.03\pm0.02$e-4 & $8.02\pm0.02$e-5 \\
	 \hline
\end{tabular}
\end{table}

For $c= 0.01,$ we have $\mathcal{R}_c(\delta_{I}(c))=0.204$, $\mathcal{R}_c(\delta_{II}(c))=0.949$, $\mathcal{R}_c(\delta_{III}(c))=0.465$.  Hence, the test $\{\delta_{II}(c)\}$ with the intuitive choice of the quantizer in  (\ref{equ:StatQntzr}) leads to a poor performance in terms of Bayes risk. Meanwhile, the test $\{\delta_{III}(c)\}$ with the ``best" invariant quantizer in (\ref{equ:BstStatQntzr}) has a better performance, and our proposed test $\{\delta_{I}(c)\}$ in Section III is the best among all three tests.

\section{Conclusion}

In this article, the problem of decentralized testing composite hypotheses in sensor networks is studied through a concrete example on testing a normal mean. Asymptotically Bayes tests $\{\delta_I(c)\}$  are constructed through a characterization of ULQ's. Contrary to our intuition, the quantizers are still of the form $I(X_{n} > \lambda).$ By exploiting the symmetries, we also investigate invariant stationary SPRT's. Numerical simulations confirm the significant advantages of our proposed test $\{\delta_I(c)\}$.

While we only consider a special problem of testing a normal mean, the essential ideas can be easily extended to other general distributions or the problem of testing $K$ ($K \ge 3$) hypotheses. It will be interesting to understand when is possible to characterize the ULQ's as in  (\ref{equ:ULQcharacter}). Another natural extension is to study the networks with multiple sensors, where different sensors may use different quantizers. The details will be investigated in our future research.

\section*{Acknowledgment}
This work was supported in part by the AFOSR grant FA9550-08-1-0376 and the NSF Grant CCF-0830472.

\end{document}